\documentclass[11pt]{article}
\usepackage{url}

\usepackage{pb-diagram,pb-xy}
\usepackage[cmtip,arrow]{xy}
\usepackage{enumerate}
\usepackage{amssymb}
\usepackage{amsmath}
\usepackage{latexsym}
\usepackage{pdfpages}

\usepackage{amscd}
\usepackage{amsfonts}
\usepackage{theorem}

\input xy
\xyoption{all}

\theoremstyle{change}
\newtheorem{Thm}{Theorem}[section]
\newtheorem{Cor}[Thm]{Corollary}
\newtheorem{Lem}[Thm]{Lemma}
\newtheorem{Pro}[Thm]{Proposition}

{\theorembodyfont{\rmfamily} \theoremstyle{change}
                             \newtheorem{Rem}[Thm]{Remark}
                             \newtheorem{Def}[Thm]{Definition}
                             
}
\newcommand{\ZZ}{\mathbb {Z}}

\newcommand{\Proof} {\noindent{\itshape Proof.\quad }}

\newcommand{\qed}{\quad\hfill$\square$}
\newcommand{\syl}{Sylow $p$-subgroup }
\newcommand{\syls}{Sylow $p$-subgroups }

\date{}

\begin{document}

\title{A GENERALIZED GOURSAT LEMMA}


\footnotetext[2] {The second named author was supported during this work by
a travel grant from the Indian Statistical Institute, Kolkata.}

\footnotetext[3] {The third named author was supported during
this work by a Discovery Grant from the Natural Sciences and
Engineering Research Council of Canada.}

\author{ K. Bauer,  D. Sen,  P. Zvengrowski }

\maketitle

\begin{abstract}


In this note the usual Goursat lemma, which describes subgroups of
the direct product of two groups, is generalized to describing subgroups
of a direct product \ $A_1\times A_2 \times \cdots \times A_n$ \ of a 
finite number of groups. Other possible generalizations are discussed and
 applications characterizing several types of subgroups are given. Most
of these applications are straightforward, while somewhat deeper applications
occur in the case of profinite groups,
 cyclic groups, and the Sylow $p$- subgroups (including
infinite groups that are virtual $p$-groups).

\end{abstract}

\section{Introduction} \label{sec:1}

In Sections 11-12 of a paper written in 1889 \cite{Goursat}, the famed French mathematician $\acute{\mbox{E}}$douard Goursat developed what is now called Goursat's lemma (also called Goursat's theorem or Goursat's other theorem), for characterizing the subgroups of the direct product $A\times B$ of two groups $A,B$. It seems to have been first attributed to Goursat by J. Lambek in \cite{Lambek2, Lambek1}, who in turn attributes H.S.M. Coxeter for bringing this to his attention. The lemma is elementary and a fundamental question to consider, for example it appears as Exercise 5, p. 75, in Lang's Algebra \cite{Lang}. It has also been the subject of recent expository articles \cite{Petrillo1, Petrillo2} in an undergraduate mathematics journal. It is possible that other authors discovered the lemma independently without knowing the original reference. Indeed, one such example, related to the theory of Lie groups,
 occurs in 1961 in a paper of A. Hattori, (\cite[Section 2.3]{Hattori}, now translated into English \cite{Z} from the original Japanese.

Other sources that mention Goursat's lemma include papers of S. Dickson \cite{Dickson} in 1969, K. Ribet \cite{Ribet} in 1976, a book by R. Schmidt \cite[Chapter 1.6]{Schmidt}, a paper of D. Anderson and V. Camillo in 2009 \cite{Anderson1}, a preprint of A. Greicius in 2009 \cite{Greicius}, a recent preprint by L. T\'oth \cite{Toth},  and several internet sites such as \cite{mathforum,poster}. Taken together, these various sources demonstrate the applicability of Goursat's
lemma to diverse branches of mathematics.

There are a number of interesting possibilities for generalizing this useful lemma. The first is to subgroups of a semi-direct product, and this is studied in \cite{Usenko}. The second is to other categories besides groups. Indeed, it is proved for modules in \cite {Lambek1}, and this implies that it will hold in any abelian category by applying the embedding theorems of Lubkin-Freyd-Heron-Mitchell cf. \cite[p. 205]{Maclane}. It is proved for rings in \cite{Anderson1}. The most general category in which one can hope to have a Goursat lemma is  an exact Mal'cev category, cf. \cite{mathforum}, and for a proof of this fact cf. \cite[Theorem 5.7]{Carboni}.

In this note we examine another generalization, to the direct product of a finite number of groups. While this seems at first glance to be a triviality since we can write $A\times B\times C\approx (A\times B)\times C$, unexpected complications arise as noted by Arroyo et al. in \cite{poster}. The complications are overcome by considering an asymmetric version of the lemma (cf. \cite{Schmidt}) in Section \ref{sec:2}, which enables us to solve the general case in Section \ref{sec:3}.   Applications within group theory are
given in Section \ref{sec:4}. These are divided into relatively easy applications followed by a
more subtle application to profinite groups, then to the  cyclic subgroups of a direct product $A_1\times \cdots \times A_n$, and finally to the Sylow $p$-subgroups of a subgroup $G$  of a direct product. The conditions for Sylow $p$-subgroups are obtained not only for finite
groups but also for virtual $p$-groups  (groups having a subgroup that is a $p$-group and has finite index). 
An Appendix gives an example that illustrates the necessity for the  
asymmetric version of Goursat's lemma in order to generalize it to finite direct products.

\noindent{\bf Acknowledgements:}    We thank L\'asl\'o T\'oth for providing us with an alternate proof of Lemma \ref{l:dirichlet} which uses the Chinese Remainder Theorem, as well as for pointing out a gap in the original proof of Theorem \ref{4.2}. We also thank Joseph Rotman and Hanafi 
Farahat for further
discussions that have been helpful with the group theoretic applications in Section 4.


\section{ Goursat's lemma, two versions} \label{sec:2}
 For convenience, and to establish the notation, we start by stating the usual (symmetric)  version of Goursat's lemma. Let $A,B$ be groups and $G\leq A\times B$ be a subgroup. The neutral element of each group $A$ and $B$, with slight abuse of notation, will be written `e'. Let $\pi_1\colon A\times B\rightarrow A,~\pi_2\colon A\times B\rightarrow B$ be the natural projections and $\imath_1\colon A\rightarrow A\times B ,~\imath_2\colon B\rightarrow A\times B$ be the usual inclusions.  
\begin{Thm}[Goursat's lemma]\label{2.1}
 There is a bijective correspondence between subgroups $G$ of $A\times B$ and quintuples $\{ \overline{G}_1, G_1, \overline{G}_2, G_2, \theta\}$, where 
$G_1\unlhd \overline{G}_1\leq A$,
 $G_2\unlhd \overline{G}_2\leq B$, and $\theta\colon \overline{G}_1/G_1\xrightarrow{\approx} \overline{G}_2/G_2$ is an isomorphism.
\end{Thm}

As mentioned above, the proof is elementary and given as an exercise in $\cite{Lang},$ it can also be found in \cite{Anderson1,Hattori}.  The basic idea of the proof is as follows.  Suppose that $G$ is a subgroup of $A\times B$.   Write $\overline{G}_1=\pi_1(G)=\{a\in A|(a,b)\in G ~\mbox{for some}~ b\in B \}$, $G_1=\imath_1^{-1}(G)=\{a\in A|(a,e)\in G\}$, similarly for $\overline{G}_2$ and $G_2$. It is easily seen that $G_1\unlhd \overline{G}_1$, $G_2\unlhd \overline{G}_2$.  The isomorphism $\theta\colon \overline{G}_1/G_1\xrightarrow{\approx} \overline{G}_2/G_2$ is  given by $\theta([a])=[b]$, where $(a,b)\in G$ and $[a]=G_1a$, $[b]=G_2b$ are the respective cosets of $a$ and $b$ in $\overline{G}_1/G_1$, $\overline{G}_2/G_2$ (again with slight abuse of notation). It is easily checked that $\theta$ is independent of the
choices of $a, b$.  Thus $G$ determines the quintuple $Q_2'(G)=\{\overline{G}_1,G_1,\overline{G}_2,G_2,\theta\}$.

Conversely, given a quintuple  $Q'=\{\overline{G}_1,G_1,\overline{G}_2,G_2,\theta\}$ where $G_1\unlhd \overline{G}_1\leq A$, $G_2\unlhd \overline{G}_2\leq B$ and $\theta\colon \overline{G}_1/G_1\xrightarrow{\approx} \overline{G}_2/G_2$,  let $\Gamma_2'(Q')$ be the subgroup $p^{-1}(\mathcal{G}_{\theta}),$ where $\mathcal{G}_{\theta}\leq \overline{G}_1/G_1\times \overline{G}_2/G_2$ is the graph of $\theta$ and $p\colon \overline{G}_1\times \overline{G}_2\rightarrow \overline{G}_1/G_1\times \overline{G}_2/G_2$ is the natural surjection.  The functions $Q_2'$ and $\Gamma_2'$ are inverse to each other.

\begin{Def} \label{quint}  Motivated by the correspondence between subgroups $G$ of $A\times B$ and the quintuples $Q'$, we say that the quintuple $Q'_2(G)$ of Theorem \ref{2.1} is the {\it Goursat quintuple} for $G$. \end{Def}

We now state an equivalent asymmetric version of the lemma cf. \cite[Theorem 1.6.1]{Schmidt}, which is in effect a minor variation of Theorem \ref{2.1} but has the advantage that it generalizes easily to higher direct products, as we shall see in Section \ref{sec:3}. An example in the Appendix shows why this asymmetric version is necessary to deal with the higher direct products.

\begin{Thm}[Asymmetric version of Goursat's lemma]\label{2.2}
There is a bijective correspondence between subgroups $G$ of $A\times B$ and quadruples $\{\overline{G}_1,\overline{G}_2,G_2,\theta_1\}$, where $\overline{G}_1\leq A$, $G_2\unlhd \overline{G}_2\leq B$ are arbitrary subgroups of $A$ and $B$, and $\theta_1\colon \overline{G}_1\twoheadrightarrow \overline{G}_2/G_2$ is a surjective homomorphism.
\end{Thm}

\Proof The proof is similar to that of Theorem 2.1.   For any subgroup $G$ of $A\times B$ we define $Q_2(G)$ to be the quadruple
\[  Q_2(G) : =  \{\overline{G}_1,\overline{G}_2,G_2,\theta_1\},\]
where $\overline{G}_1$, $\overline{G}_2$ and $G_2$ are the first, third and fourth coordinates of $Q_2'(G)$ (cf. Theorem 2.1).  The surjection $\theta_1$ is given by $\theta_1(a)=[b]$ for $a\in \overline{G}_1$ and $(a,b)\in G$ for some $b\in \overline{G}_2$ (again easily seen to be independent of the choice of $b$).  

Conversely, for an arbitrary quadruple $Q=\{\overline{G}_1,\overline{G}_2,G_2,\theta_1\}$, with $\overline{G}_1\leq A$, $G_2\unlhd \overline{G}_2\leq B$, and $\theta_1\colon \overline{G}_1\twoheadrightarrow \overline{G}_2/G_2$ a surjective homomorphism define 
\[\Gamma_2(Q):=p^{-1}(\mathcal{G}_{\theta_1}),\]
 where $\mathcal{G}_{\theta_1}\subseteq \overline{G}_1\times \left(\overline{G}_2/G_2\right)$ is the graph of $\theta_1$ and $p\colon \overline{G}_1\times \overline{G}_2\rightarrow \overline{G}_1\times \left(\overline{G}_2/G_2\right)$ is the natural surjection.  The functions $Q_2$ and $\Gamma_2$ are inverse to each other.
 
 \qed

The equivalence of Theorem \ref{2.1} and Theorem \ref{2.2} is easily seen. It need only be pointed out that $\theta$ determines the surjection $\theta_1$ as the composition $$\overline{G}_1\twoheadrightarrow \overline{G}_1/G_1\xrightarrow[\approx]{\theta} \overline{G}_2/G_2,$$ while $\theta_1$ determines $\theta$ via the first isomorphism theorem, specifically
\[
 \begin{diagram}
\node{\overline{G}_1}\arrow{e,t}{\theta_1}\arrow{s,A}\node{\overline{G}_2/G_2}\\
\node{\overline{G}_1/\mathrm{Ker}(\theta_1)}\arrow{ne,tb,A}{\approx}{\theta}
 \end{diagram}
.
\]
Finally, we have
\begin{equation*}
 \begin{split}
\mathrm{Ker}( \theta_1)&=\{a\in \overline{G}_1|\theta_1(a)=[b], \mbox{with}~b\in G_2\}\\
& =\{a\in \overline{G}_1|(a,b)\in G,~(e,b)\in G \}\\
&= \{a\in \overline{G}_1|(a,e)\in G\}\\
&= G_1.
 \end{split}
\end{equation*}

\begin{Rem}\label {2.3}
It is often useful to note a few additional facts, which we now list. For (a),
(b) cf. \cite{Hattori}, \cite{Z}.

(a) \ \  $(G_1\times G_2)\unlhd G\leq \overline{G}_1\times\overline{G}_2$,

(b) \ \  $G/(G_1\times G_2)\approx \overline{G}_1/G_1\approx  
\overline{G}_2/G_2$,  

(c) \ \  From (b),  using \ $|G_1\times G_2| = |G_1|\cdot |G_2|$, one readily obtains that
 \ $|G| = |G_1|\cdot |\overline{G}_2|
 =  |G_2|\cdot |\overline{G}_1|$ \ for the cardinality of $G$,

(d) \ \  One has short exact sequences of groups

\[
 \begin{diagram}
\node{\iota_2(G_2)= {\rm Ker}(\pi_1|G)}\arrow[2]{e,t,J}{\unlhd}
\node[2]{G}\arrow[2]{e,t,A}{\pi_1 |G}\node[2]{\overline{G}_1 \ ,}    \\
\mbox{} \\
\node{\iota_1(G_1)= {\rm Ker}(\pi_2|G)}\arrow[2]{e,t,J}{\unlhd}
\node[2]{G}\arrow[2]{e,t,A}{\pi_2 |G}\node[2]{\overline{G}_2 \ .}
 \end{diagram}
\] 

\end{Rem}


\section{The generalized Goursat lemma } \label{sec:3}

As mentioned in the Introduction, generalizing the usual Goursat lemma from $n=2$ to $n\geq 2$ seems to create unexpected complications. However, using the asymmetric version Theorem $2.2$, the generalization to finite $n\geq 2$ becomes routine. We will state the result (Theorem 3.2 below) for $n\geq 2$,
after
first introducing some convenient notation for any subgroup $G$ of  $A_1\times \cdots \times A_n$.

\begin{Def}\label{d:3.1}
 Let $S\subset \{1,2,\cdots,n\}=\underline{n}$, and $j\in \underline{n}\smallsetminus S$. Then $$G(j|S):=\{x_j\in A_j|(x_1,\cdots,x_j,\cdots,x_n)\in G~~~~~~~~~~~~~~~~~~~~~~~~~~~~~~~~~~~~~~~~$$ $$~~~~~\mbox{for some}~x_i\in A_i, 1\leq i \leq n, i\neq j,
 ~\mbox{with}~x_i=e ~\mbox{if}~i\in S \}.$$
\end{Def}
For example, $$G(1|\emptyset)=\pi_1(G),~~G(1|\{2,3,\cdots, n-1\})=\{x_1\in A_1|(x_1,e,\cdots,e)\in G\}.$$ These correspond to the notation used in Section $2$, when $n=2$, via $G(1|\emptyset)=\overline{G}_1, ~G(1|\{2\})=G_1, ~G(2|\emptyset)=\overline{G}_2, ~G(2|\{1\})=G_2.$   For brevity, we extend this notation and let $\overline{G}_k :~= G(k|\emptyset)$ for all $k$.  
For convenience, we shall usually omit the brackets $\{\}$, e.g. $G(1|\{2,3\})=G(1|2,3)$. Note that if $T\subseteq S$, then $G(j|S)\lhd G(j|T)$.
As in Section \ref{sec:2} we let $\overline{G}_i=\pi_i(G)$, where $\pi_i:A_1\times \cdots \times A_n \twoheadrightarrow A_i$ is the standard projection onto the $i$-th factor.
Finally, it will be convenient to also use \ $\Pi_i : A_1 \times \cdots \times A_n 
\twoheadrightarrow  A_1 \times \cdots \times A_i$ \ for the standard projection onto the
first $i$ factors, \ $1 \leq i \leq n$ \ (e.g. \ $\Pi_1 = \pi_1$ \ and \ $\Pi_n = 
{\rm id}_{A_1 \times \cdots \times A_n} $).

\begin{Thm}[Goursat's lemma for $n \geq 2$]\label{3.2}
There is a bijective correspondence between the subgroups \ $G \leq A_1 \times \cdots\times A_n$ and $(3n-2)$-tuples 
$$Q_n(G) :=    \{\overline{G}_1,\overline{G}_2, G(2|1), \theta_1,
\ldots, \overline{G}_n, G(n|1,\ldots ,n-1), \theta_{n-1}  \} , $$
where $\overline{G}_i \leq A_i, \ G(i|1,\ldots,i-1)\unlhd \overline{G}_i,   $ and
\ $\theta_i : \Lambda_i \twoheadrightarrow \overline{G}_{i+1}/G(i+1|1,\ldots,i)$ is a surjective
homomorphism. Here
$\Lambda_i \leq A_1\times \cdots\times A_i$ is defined recursively, $1 \leq i \leq n-1$, by setting \ $\Lambda_1  := \overline{G}_1$ \ and 
 $$\Lambda_{i+1} := \Gamma_2(\{ \Lambda_i, \overline{G}_{i+1}, G(i+1|1,\ldots ,i), \theta_i \}) \leq
 (A_1\times \cdots \times A_i)\times A_{i+1} \ ,$$
 with  $\Gamma_2$ as defined in Theorem 2.2.  
 \end{Thm}
 
\Proof
Starting with $G$, we must construct the $(3n-2)$-tuple $Q_{n}(G)$, all entries of which are already
defined (from $G$) except the $\theta_i$. We shall show by induction, $1\leq i \leq n-1$, first that $\Lambda_i
 = \Pi_i(G)$, and second that $\theta_i$ can then be suitably defined to successfully carry out
 the inductive step.
 
 To start the induction we simply observe that by hypothesis \  $\Lambda_1 = \overline{G}_1$,
 \ and \ $\overline{G}_1 = \pi_1(G) = \Pi_1(G)$. Now suppose, as inductive
hypothesis, that \ $\Lambda_i = \Pi_i(G)$. To define \ $\theta_i : \Lambda_i \twoheadrightarrow \overline{G}_{i+1}/G(i+1|1,\ldots,i)$, suppose \ $x \in \Lambda_i$. By the
inductive hypothesis \ $x = (a_1,\ldots ,a_i) \in 
\Pi_i(G) \leq A_1\times \cdots \times A_i$. Then \ $(a_1,\ldots a_i, a_{i+1}, 
\ldots, a_n) \in G$ \ for some \ $a_j \in A_j, \ j = i+1,\ldots, n$. We define
\ $\theta_i(x) := [a_{i+1}] \in \overline{G}_{i+1}/G(i+1|1,\ldots,i)$. To see
that this definition makes sense one must check that \ $a_{i+1}\in 
\overline{G}_{i+1}$, \ that $\theta_i$ is surjective, and that the definition is
independent of the choice of $a_{i+1}$. The first two are obvious, and as far as
the third suppose \ $(a_1,\ldots,a_i,a'_i,\ldots,a'_n)\in G$. Then \ 
$(e,\ldots,e,a_{i+1}^{-1}a'_{i+1},\ldots,a_n^{-1}a'_n) \in G$ \ which implies
\ $a_{i+1}^{-1}a'_{i+1} \in  G(i+1|1,\ldots,i)$ and therefore \ $[a_{i+1}]
 = [a'_{i+1}]$, \ i.e. $\theta_i$ is well defined. From the definitions of 
$\Lambda_{i+1}$, 
of $\Gamma_2$  (Theorem 2.2), and of $\theta_i$,  it 
is now obvious that \  $\Lambda_{i+1} = \Pi_{i+1}(G)$. This completes the 
inductive step.

Conversely, suppose we are given a $(3n-2)$-tuple
 $$Q =    \{\overline{G}_1,\overline{G}_2, G(2|1), \theta_1,\overline{G}_3, G(3|1,2), \theta_2 ,
\ldots, \overline{G}_n, G(n|1,\ldots ,n-1), \theta_{n-1}  \} , $$
satisfying the hypotheses of the theorem. Then simply define a subgroup \ $\Gamma_n(Q) \leq 
A_1\times \cdots\times A_n$ \ by \ $\Gamma_n(Q) := \Lambda_n$. Because $Q_2$ and $\Gamma_2$ are
inverse to one another (cf. proof of Theorem 2.2), the iterated versions of these two
operations, namely $Q_n$ and $\Gamma_n$, are also inverse to one another. \hfill $\Box$ 

 \begin{Def} For a subgroup $G\leq A_1\times \cdots \times A_n$, we say that the corresponding $(3n-2)$-tuple $Q_n(G)$ of Theorem 3.2 is the {\it Goursat decomposition} of $G$. \end{Def}

We may also refer to the Goursat quintuple of Definition \ref{quint} as a Goursat decomposition, even though it differs slightly from the quadruple $Q_2(G)$.  In this paper the context makes clear the difference between these two possible decompositions.

\begin{Rem}\label {3.3}
As in Remark 2.3, \ $G \leq \overline{G}_1\times \cdots \times \overline{G}_n$. However the rest
of Remark 2.3 does not directly apply here.

\end{Rem}

\mbox{}

\bigskip


\section{Applications} \label{sec:4}
There are many potentially interesting applications of the (generalized) Goursat lemma within group theory. In this section we start with several
 easy applications and then explore three
relatively deeper applications :  profinite groups, 
cyclic groups and $p$-Sylow subgroups. For example, an immediate consequence of the lemma is that the subgroup $G$ is a sub-direct product if and only if $\theta_j$ is the trivial homomorphism, \ $1 \leq j \leq n-1$.

\begin{Pro}\label{4.1}
Let ${\cal C}$ be a class of groups closed under taking subgroups, quotient groups, and finite
direct products. Let $G$ be a subgroup of $A_1\times\cdots\times A_n$. Then $G$ is in the class ${\cal C}$ iff each
$\overline{G}_i$ is in ${\cal C}$.

\end{Pro}

\Proof The proposition follows from Theorem \ref{3.2} and Remark \ref{3.3}, 
which taken together indicate that $G$ is a subgroup of 
$ \overline{G}_1 \times \ldots \times \overline{G}_n$ and furthermore each $\overline{G}_i$ is a quotient of $G$.
\qed

For example, Proposition \ref{4.1} holds for any Serre class of abelian groups and
any variety of groups (in the sense of \cite{Neumann}). In particular it holds
for each of the following familiar families of groups : \ (a) \ finite groups,
\ (b) \ abelian groups, \ (c) \ $p$-groups, \ (d) \ nilpotent groups, 
 \ (e) \ groups nilpotent of class at most $m$, \ (e) \ solvable groups.

%
To verify each of these examples, one need only verify that each class of groups satisfies the hypotheses of Proposition \ref{4.1}. This is trivial
for (a), (b), (c).   In (e) we find that groups which are nilpotent of class at most $m$ are closed under subgroups and quotients as \cite[Theorem 5.35, 5.36]{Rotman}.  As well,  the products of groups which are nilpotent of class at most $m$ must be nilpotent of class at most $m$ since the commutator subgroup of the product is the product of the commutator subgroups (similarly for (d)). For (f)
cf. \cite[Theorems 6.11, 6.12, and Corollary 6.14]{Rotman}.


It is interesting to note that Hattori, in \cite{Hattori} or \cite{Z}, first determined the finite subgroups of the Lie group $S^3$ (which is isomorphic to SU($2$)$\approx$ Sp($1$)$\approx$ Spin($3$) as a Lie group). He then applied Goursat's lemma and Proposition \ref{4.1} to determine all finite subgroups of $S^3\times S^3$. Using the results of Section \ref{sec:3} we could now, for example, find all finite subgroups of $S^3\times \underset{n}{\cdots}\times S^3$.  It is also interesting that, in fact, the  papers of Goursat \cite{Goursat} and Hattori \cite{Hattori} or \cite{Z} study closely related questions.

As another interesting application, which also involves some topology, we consider profinite groups.  Their
 definition and basic properties can be
 found in \cite{Ramakrishnan}.  Briefly, a topological group is a profinite group if it can be obtained as an inverse limit of finite groups, each having the discrete topology. Profinite groups can also be characterized as topological
groups that are Hausdorff, compact, and totally disconnected.  The class of profinite groups is closed with respect to taking closed subgroups, quotient groups by a closed normal subgroup, and (arbitrary) direct products.  However, arbitrary subgroups of profinite groups may not be profinite, so Proposition \ref{4.1} does not immediately apply.  We shall nevertheless be able obtain a similar result by being careful about the topology. 

For simplicity we start with two profinite groups $A, B$ and a subgroup $G \leq A \times B$, where $A\times B$ has the product topology and $G$ the subspace topology. Of course $G$ has a Goursat decomposition \ $Q_2(G) = \{ \overline{G}_1,  \overline{G}_2, G_2, \theta_1 \}$ \ as a group.
We topologize $\overline{G}_i, i=1,2$, using the surjection \ $\pi_i : G \twoheadrightarrow 
\overline{G}_i$ and giving $\overline{G}_i$ the identification topology (also called the quotient topology).
 Then $G_1 \unlhd \overline{G}_1 $ is given the subspace topology. We also note that
the usual projection and inclusion maps $\pi_i, \ \iota_i$ are continuous, since the product topology is being used, and that each space $\overline{G}_i, G_i$ is a subspace of either $A$ or
$B$, hence is Hausdorff. Note that the notation $\overline{G}_i$ is being used here as in
Theorem 2.1 (definition) and subsequently, it has nothing to do with the closure operator in
topology.

\begin{Pro}  Let $A, B$ be profinite groups and $G \leq A\times B$ as above.  Then $G$ is a profinite group iff each of the subgroups  in the Goursat decomposition for $G$ are profinite groups and $\theta_1$ is continuous.  

\end{Pro}

\Proof  Suppose that $G$ is a profinite group. Using the Hausdorff property we have
$(\{ e\}\times B)\bigcap G$ is a closed subgroup of $G$, so also profinite. But then \ $\pi_2 : 
  (\{ e\}\times B)\bigcap G  \to G_2  $ \ is a continuous bijection of a compact space onto
   a Hausdorff space,
hence a homeomorphism. Since it is also a group isomorphism, $G_2$ is profinite, and similarly for $G_1.$
 Since $\iota_i(G_i)$ are closed normal subgroups of $G$, Remark 2.4(d)
 and the properties
of profinite groups imply that $\overline{G}_i$ are also profinite. Now consider the
continuous surjective homomorphism \ 
${\rm id}\times p : \overline{G}_1\times \overline{G}_2 \twoheadrightarrow 
\overline{G}_1\times (\overline{G}_2/G_2).$ Recalling that \ $G = 
({\rm id}\times p)^{-1}({\cal G}_{\theta_1})$, we have  ${\cal G}
(\theta_1)= ({\rm id}\times p) (G)$. This is a continuous image of a compact
space, hence compact, and hence closed since it lies in the Hausdorff space
$\overline{G}_1\times (\overline{G}_2/G_2).$ By the closed graph theorem
$\theta_1$ is continuous.

Conversely, suppose that $\overline{G}_1$, $\overline{G}_2$, $G_2$ are profinite groups, and $\theta_1$ is continuous.  Then $\overline{G}_2/G_2$ is also a profinite group.  Since $\theta_1$ is continuous, we must have (again by the 
closed graph theorem) that the graph ${\cal G}_{\theta_1} \subseteq
\overline{G}_1\times (\overline{G}_2/G_2)$
 is a closed subgroup.  The identification
 map $q:
  \overline{G}_1\times \overline{G}_2\to \overline{G}_1\times (\overline{G}_2/G_2)$ is continuous,  therefore $G=q^{-1}({\cal G}_{\theta_1})$
   is a closed subgroup of the profinite group $\overline{G}_1\times \overline{G}_2$.  Hence $G$, being a closed subgroup of a profinite group, is itself profinite. 
\qed

We remark that there is a close relation between Proposition 4.2 above and 
 \cite[Lemma 4.6]{Greicius}, although neither one implies the other. It is also clear that
 Proposition 4.2 will generalize to subgroups of finite direct products of profinite groups with more than two factors,
 in the obvious way.

 The next application, that of determining the cyclic subgroups of $A\times B$, will involve more substantial use of Goursat's lemma.  Cyclic subgroups are not closed under products, so Proposition \ref{4.1} does not apply.   We shall henceforth use additive notation since $\overline{G}_1,G_1,\overline{G}_2,G_2$ will be abelian.  One preliminary lemma will be needed.
 
 \begin{Lem}\label{l:dirichlet}  Let $d$ divide both $m$ and $n$, and let
 \[ \theta:{\mathbb Z}_m/d{\mathbb Z}_m\to {\mathbb Z}_n/d{\mathbb Z}_n\]
 be a given isomorphism, where both groups are of course isomorphic to ${\mathbb Z}_d$.  Then there exist generators $\alpha$, $\beta$ of respectively ${\mathbb Z}_m$ and ${\mathbb Z}_n$ such that $\theta([\alpha])=[\beta]$, where $[\alpha]=\alpha+d{\mathbb Z}_m$ and similarly $[\beta]=\beta+d{\mathbb Z}_n$.
 \end{Lem}
 
 \Proof We take $\alpha=1$, and set $\theta([1])=[\beta_0]$, so $[\beta_0]$ generates ${\mathbb Z}_n/d{\mathbb Z}_n\cong {\mathbb Z}_d$.  This means that $d$ and $\beta_0$ are coprime.  In general $\beta_0$ and $n$ need not be coprime, however $[\beta_0]=[\beta_0+d]=[\beta_0+2d]= \cdots$, so it will suffice to show that $\beta_0+kd$ is coprime to $n$ for some $k$.  By Dirichlet's famous theorem \cite{Dirichlet} there are infinitely many primes in the arithmetic progression $\{\beta_0+kd\}$, so we can choose $k$ with $p=\beta_0+kd$ prime and also $p$ not a divisor of $n$.  Then $p$ and $n$ are coprime, and taking $\beta=p\in {\mathbb Z}_n$ fulfills the conclusion of the lemma.
 
 \qed
 
Lemma \ref{l:dirichlet}  is necessary because given  cyclic groups $G_1$ and $G_2$, subgroups $H_1 \unlhd G_1$ and $H_2 \unlhd G_2$, and an isomorphism $\theta: G_1/H_1\to G_2/H_2$,  
a representative $\beta \in G_2$ of
$\theta([\alpha])$, $\alpha$ being a generator of $G_1$, need not generate $G_2$.  
For example, there is an isomorphism
 \[ \theta: {\mathbb Z}_{45}/9{\mathbb Z}_{45} \to {\mathbb {Z}}_{198}/9{\mathbb Z}_{198}\]
 defined by $\theta([1])=[2]$.  But $2\in {\mathbb {Z}}_{198}$ is not a generator.  However, $[2]=[11]=[20]=[29]$ and although $11$ and $20$ are also not generators of ${\mathbb Z}_{198}$, 29 is coprime to 198 and so is a generator.

\begin{Thm}\label{4.2}
 Let $G$ be a subgroup of $A\times B$ with  Goursat quintuple $Q'_2(G)=\{ \overline{G}_1, G_1, \overline{G}_2, G_2, \theta\}$ .
\begin{enumerate}[(a)]
 \item The subgroup $G$ is finite cyclic if and only if \  $\overline{G}_1, \overline{G}_2$ are both finite cyclic and $G_1,G_2$ have coprime order. Furthermore, $|G|=\mbox{lcm}(|\overline{G}_1|,|\overline{G}_2|)$.
\item The subgroup $G$ is infinite cyclic if and only if either $\overline{G}_1$ is infinite cyclic, $\overline{G}_2$ is finite cyclic, and $G_2=\{0\}$, or $\overline{G}_2$ is infinite cyclic, $\overline{G}_1$ is finite cyclic, and $G_1=\{0\}$, or both $\overline{G}_1,\overline{G}_2$ are infinite cyclic with $G_1=G_2=\{0\}$.
 \end{enumerate}
\end{Thm}
\Proof
In either case $(a)$ or $(b)$, $\overline{G}_1$ or $\overline{G}_2$ being the homomorphic image of a cyclic group, must be cyclic, hence their respective subgroups $G_1,G_2$ are also cyclic. 

(a)  Suppose $G$ is finite cyclic, then it is generated by an element $(\alpha,\beta)$, whence $\overline{G}_1$ is cyclic and generated by $\alpha$, $\overline{G}_2$ cyclic and generated by $\beta$. Let the respective orders of $\overline{G}_1,\overline{G}_2$ (i.e. of $\alpha,\beta$) be $m,n$, and set $d=\mbox{gcd}(m,n)$. Also write $m=m_1d,~n=n_1d$. Then $m_1,n_1$ are coprime, and there exist $x,y\in \mathbb{Z}$ with $xm_1+yn_1=1$, or equivalently $xm+yn=d$. Now $n(\alpha,\beta)=(n\alpha,0)\in G$ implies $n\alpha\in G_1$. Also $m\alpha=0\in G_1$. Hence $d\alpha=(xm+yn)\alpha=x(m\alpha)+y(n\alpha)\in G_1$. It follows that $c\alpha\notin G_1$ if $0<c<d$. For, if $c\alpha\in G_1$, then $(c\alpha,0)\in G$. Hence $(c\alpha,0)=z(\alpha,\beta)$ for some integer $z$. Therefore $z\beta=0=(c\mbox{-}z)\alpha$, whence $n|z$ and $m|(c\mbox{-}z)$. Since $d$ divides $m$ and $n$, we have $d$ divides both $z,c\mbox{-}z$. As a result, $d|c$, which is a contradiction. Thus we conclude that $G_1$ is the cyclic subgroup of $\overline{G}_1$, generated by $d\alpha$ and having order $m/d=m_1$. Similarly, $G_2$ is generated by $d\beta$ and has order $n_1$, so the orders of $G_1$ and $G_2$ are coprime.

Conversely, suppose $|G_1|=m_1$ is coprime to $|G_2|=n_1$, and set $d=|\overline{G}_1|/|G_1|=|\overline{G}_2|/|G_2|$, $m=m_1d,n=n_1d$. Then $d=\mbox{gcd}(m,n)$. Also $\overline{G}_1$ will be cyclic of order $m$, $\overline{G}_2$ cyclic of order $n$. Using the isomorphism $\theta\colon \overline{G}_1/G_1\rightarrow \overline{G}_2/G_2$ and Lemma \ref{l:dirichlet}, choose  generators $\alpha$ of $\overline{G}_1$ and $\beta$ of  $\overline{G}_2$ with $\theta([\alpha])=[\beta]$. Then  $[\alpha]$ and $[\beta]$ are elements of of order $d$, whence $\beta$ has order $d|G_2|=dn_1=n$ and generates $\overline{G}_2$. Also, by Goursat's lemma, $\gamma=(\alpha,\beta)\in G$. Finally, the order of $\gamma$ is lcm($o(\alpha),o(\beta)$)=lcm($m,n)=(mn/d)=dm_1n_1$. Further, again using Goursat's lemma, $|G|=|\mathcal{G}_{\theta}||G_1||G_2|=dm_1n_1$. Thus $G$ is cyclic of this order and generated by $(\alpha,\beta)$.

(b) Since $G$ is infinite and $G\subseteq \overline{G}_1\times \overline{G}_2$, at least one of $\overline{G}_1,\overline{G}_2$ must be infinite cyclic. Without loss of generality, suppose $\overline{G}_1\approx \mathbb{Z}$. Now suppose $(\alpha,\beta)$ generates the cyclic group $G$, then $\alpha$ generates $\overline{G}_1$, and $\beta$ generates $\overline{G}_2$. We claim that $G_2=\{0\}$. For, if $y\in G_2$ then $y=r\beta$ for some integer $r$, whence $(0,y)=r(0,\beta)\in G$. This implies $(0,y)=k(\alpha,\beta)=(k\alpha,k\beta)$ for some integer $k$. Therefore $k\alpha=0$, whence $k=0$ and $y=k\beta=0$. Hence $G_2=\{0\}.$ We now consider separately the cases $\overline{G}_2$ infinite and $\overline{G}_2$ finite (the case $\overline{G}_1$ finite and $\overline{G}_2\approx \mathbb{Z}$ is symmetric to the latter, so can be omitted).

Suppose first $G\approx \mathbb{Z}$ with $\overline{G}_2\approx \mathbb{Z}$. Then the argument in the previous paragraph now also implies $G_1=\{0\}$. Conversely, suppose $\overline{G}_1\approx \overline{G}_2\approx \mathbb{Z}$ and $G_1=G_2=\{0\}$. Then the isomorphisms (cf. Remark 2.3) 
$$G/(G_1\times G_2)\approx \overline{G}_1/G_1\xrightarrow{\theta} \overline{G}_2/G_2$$ reduce to $G\approx \overline{G}_1\approx \overline{G}_2\approx \mathbb{Z}$.

Secondly, for the remaining case, suppose $G\approx \mathbb{Z}, \overline{G}_1\approx \mathbb{Z}$ as before and now $\overline{G}_2\approx \mathbb{Z}_n$ is cyclic of order $n$, $n\geq 2$. Then $n(\alpha,\beta)=(n\alpha,0)$ implies $n\alpha\in G_1$ and clearly $i\alpha\notin G_1$ if $i<n$. Thus $G_1\approx n\mathbb{Z}$, and as before $G_2=\{0\}$.

Conversely, suppose $\overline{G}_1\approx \mathbb{Z},G_1\approx n\mathbb{Z},\overline{G}_2\approx \mathbb{Z}_n$ and $G_2=\{0\}$. In this case we have the isomorphism $\theta\colon \overline{G}_1/G_1\rightarrow \overline{G}_2/G_2\approx \mathbb{Z}_n$. Let $\alpha\in \overline{G}_1$ with $[\alpha]$ generating $\overline{G}_1/G_1$. Then $\theta([\alpha])=[\beta]=\beta$ generates $\overline{G}_2/G_2=\overline{G}_2\approx \mathbb{Z}_n$. We claim that $G$ is generated by the single element $(\alpha,\beta)$, and thus is infinite cyclic. To see this, let $(x,y)\in G\subseteq \overline{G}_1\times \overline{G}_2,$ so $x=j\alpha,y=k\beta$ for some integers $j,k$. Furthermore $(x,y)\in G$ implies $\theta([x])=[y]=y$, which gives $k\beta=y=\theta([j\alpha])=j(\theta([\alpha]))=j\beta$. Then $j\equiv k(\mbox{mod}~n)$, so $j(\alpha,\beta)=(j\alpha,j\beta)=(j\alpha,k\beta)=(x,y)$.

\qed

\begin{Thm}
 Let $G$ be a subgroup of $A\times B\times C$  with its Goursat decomposition $\{ \overline{G}_1, \overline{G}_2, G(2|1), \theta_1, \overline{G}_3, G(3|1,2), \theta_2\}$ .
\begin{enumerate}[(a)]
 \item The subgroup $G$ is finite cyclic if and only if $\overline{G}_1,\overline{G}_2,\overline{G}_3$ are finite cyclic and each of the pairs of integers $(|G(1|2)|,|G(2|1)|),$ $(|G(1|3)|,|G(3|1)|),$ $(|G(2|3)|,|G(3|2)|)$ is coprime. In this case one also has $$|G|=\mbox{lcm}(|\overline{G}_1|,|\overline{G}_2|,|\overline{G}_3|).$$
\item The subgroup $G$ is infinite cyclic if and only if one of the following three cases (up to obvious permutation of indices) occur :
\begin{enumerate}[(i)]
\item  $\overline{G}_1\approx \mathbb{Z}$, $\overline{G}_2$ and $\overline{G}_3$ are finite cyclic, $G(2|1)=G(3|1)=\{0\}$, and $G(2|3),G(3|2)$ are coprime.
 \item $\overline{G}_1\approx \overline{G}_2\approx \mathbb{Z}$, $\overline{G}_3$ finite cyclic, and $G(2|1)=G(3|1)=G(1|2)=G(3|2)=\{0\}$.
\item $\overline{G}_i\approx \mathbb{Z}$ for $i=1,2,3$ and $G(i|j)=0$ for $1\leq i\neq j \leq 3$.
\end{enumerate}
\end{enumerate}
 \end{Thm}
\Proof
(a) If $G$ is finite cyclic, then so is $\overline{G}_{12}=\Pi_2(G)
\subseteq A\times B$. Applying part $(a)$ of Theorem \ref{4.2} to $\overline{G}_{12}$ gives us $\overline{G}_1, \overline{G}_2$ finite cyclic with $|G(1|2)|$ coprime to $|G(2|1)|$. The other conditions follow by symmetry.

Conversely, suppose $\overline{G}_1,\overline{G}_2,\overline{G}_3$ are all 
finite cyclic with respective orders $m,n,p,$ and that the three coprimality 
conditions hold. Let \[G_{12}=\iota^{-1}_{12}(G)=\{ (a,b)\in \overline{G}_1\times 
\overline{G}_2 | (a,b,e)\in G\ \}.\]  
Applying Theorem \ref{4.2}(a) two times we obtain that $\overline{G}_{12}$ and $G_{12}$ 
are both finite cyclic with respective orders $\mbox{lcm}(m,n)$, 
$\mbox{lcm}(|G(1|3)|,|G(2|3)|)$.  
We next  apply Theorem \ref{3.2}, which tells us that $G$ is determined by the 
surjection $\theta_2\colon \overline{G}_{12}\twoheadrightarrow \overline{G}_3/G(3|1,2)$.
A third application of Theorem \ref{4.2}(a) now tells us that $G$ will be cyclic 
if $|G_{12}|$ and $|G(3|1,2)|$ are coprime. But $|G_{12}|=\mbox{lcm}(|G(1|3)|,|G(2|3)|)$,  
and $|G(3|1,2)|$ is a divisor of $G(3|1)$ which is coprime to $|G(1|3)|$. 
Thus $|G(3|1,2)|$ is coprime to $|G(1|3)|$, and similarly is coprime to $|G(2|3)|$, 
so also coprime to their least common multiple $|G_{12}|$.

(b) The three cases when $G$ is infinite cyclic all follow from Theorem \ref{4.2}(b) in obvious ways, namely in $(b)(i)$ we use $A\times B\times C\approx A\times (B\times C)$, in $(b)(ii)$ and $(b)(iii)$ we use $A\times B\times C\approx (A\times B)\times C$. We omit the details.

\qed

The generalization of this theorem to $n\geq 3$ is now clear, albeit cumbersome to state since
there will be many cases involved.

\smallskip

Determining the \syls  of a a group \ $G\ \leq A_1\times\cdots \times A_n$ in terms of the Goursat
decomposition of $G$ is an application of a slightly different type. Our main result in this 
direction, Theorem 4.8 below, gives a very simple and natural answer to this question, and not only for
finite groups but for certain classes of infinite groups. We therefore commence with a brief
discussion of \syls for groups that are not necessarily finite, taking $p$ to be a fixed prime
for the remainder of this section.

The \syls of an arbitrary group $G$ are easily defined as its maximal $p$-subgroups, which always exist by a Zorn's lemma argument. However,
following \cite[Section 14.3]{Robinson}, one sees that without some sort of
finiteness hypothesis the familiar Sylow theorems (for a finite group) can fail
badly. Indeed it is possible for two \syls to even have different cardinalities, let alone be isomorphic or conjugate. One hypothesis that will insure the
usual Sylow theorems hold, namely that all \syls are conjugate and their number is
both finite and congruent to $1$ modulo $p$, is that there exists a \syl with
a finite number of conjugates. This theorem was proved in 1938 by
Dieman, Kurosh, Uztov \cite{Dieman}  and in 1940 by Baer \cite{Baer} . We shall 
call this finite conjugacy
property ``${\rm FC}_p$'' and also call this theorem the ``${\rm FC}_p$  theorem.''  We shall consider a further finiteness property, that the group
is a virtual $p$-group, i.e. it has a $p$-subgroup with finite index.     

  It is easy to see that a virtual $p$-group satisfies ${\rm FC}_p$, for suppose $G$ is a virtual
  $p$-group. Then
it has a $p$-subgroup $M$ of finite index, and without loss of generality (enlarging $M$ if necessary) we
may suppose $M$ is a \syl. The number of conjugates of $M$ is given by the index $[G:N_M]$ of its
normalizer $N_M$, and since 
$N_M \geq M$ this index is finite, proving the ${\rm FC}_p$ property. On the other 
hand, since any abelian group
trivially satisfies ${\rm FC}_p$, it is clear that ${\rm FC}_p$ does not imply that the group is
a virtual $p$-group.

The next proposition gives us an easy way to identify a \syl \ in a virtual $p$-group and along with the lemma that follows it  will make the proof
of the main result,  Theorem 4.8 below, quite easy.

\begin{Pro} Let $G$ be a virtual $p$-group and $M$ a $p$-subgroup. Then $M$ is a \syl \ if and only if
$[G:M]$ is coprime to $p$.
\end{Pro}
\begin{Proof} Suppose $M \leq H \leq G$, then \ $[G:M] = [G:H]\cdot [H:M]$. If $[G:M]$ is coprime to $p$
then it follows that $[H:M]$ is also coprime to $p$. If also $H$ is a $p$-group this can only
happen if $[H:M] = 1$, whence $M$ is a maximal $p$-subgroup and thus Sylow. 

Conversely suppose
$M_1$ is a \syl. We have already argued above that any virtual $p$-group $G$ admits a \syl \ $M$
having finite index. Since $G$ also satisfies ${\rm FC}_p, \ M, M_1$ are conjugate by the
${\rm FC}_p$ theorem and hence have the same
index in $G$, so we can deal with $M$. Now \ $[G:M] = [G:N_M] \cdot [N_M:M]$, all being finite numbers. Again by the ${\rm FC}_p$ theorem, $[G:N_M] \equiv 1 \ 
({\rm mod} \ p), $ and is thus coprime to $p$. Since $M \trianglelefteq N_M$ and is also a \syl \ of $N_M$, it is easy to see that $[N_M:M]$ is also coprime
to $p$ by considering the finite group $N_M/M : = A$  which has order $[N_M:M]$, and the surjection
\ $ N_M \twoheadrightarrow A$ (the \syl \ of $A$ must be trivial, otherwise
its inverse image would be a $p$-subgroup of $N_M$ strictly larger than $M$).
\qed           
\end{Proof}

The next lemma is a ``non-commutative'' version of a result that is 
familiar in abelian categories, where it follows at once by
taking the quotient objects to form a third exact row in the diagram.

\begin{Lem}  Let $G$ be any group, $M$ any subgroup, and $\pi$ a surjective
homomorphism of $G$ onto a group $H$. Consider the commutative diagram

\[
 \begin{diagram}
\node{K_1}\arrow{s,1,J}\arrow[2]{e,t,J}{}\node[2]{M}\arrow{s,1,J}\arrow[2]{e,t,A}{\pi |M}\node[2]{H_1}\arrow{s,1,J}   \\
\node{K}\arrow[2]{e,t,J}{} \node[2]{G} \arrow[2]{e,t,A}{\pi} \node[2]{H \ ,}  
 \end{diagram}
\]
where $K, K_1$ are the respective kernels of $\pi, \pi|M$, and $H_1$ is the image of
$\pi|M$. Then one has the following relation of (possibly infinite) cardinal numbers : \ \  
$ [G:M] = [H:H_1] \cdot [K:K_1] \ . $
\end{Lem}
\begin{Proof}
Since $K \unlhd G$, it is standard, cf. \cite[Theorem 2.13]{Rotman}, that 
$KM = MK = K\vee M $ is a subgroup of $G$. 
Another standard fact is that \ $[G:M] = [G:MK]\cdot [MK:M]$, cf. 
 \cite[Section 1.3.5]{Robinson}.
The plan of this proof is to show that \ $[G:MK] = [H:H_1]$ \ and \ $[MK:M]
 = [K:K_1]$, \ which will complete the proof. 

For \ $[G:MK]$ \ let us first note that\ $MK = \pi^{-1}
(H_1)$.  Then  $[G:MK] = [H:H_1]$ is part of
  \cite[Theorem 2.3.4]{Hall}.

For \ $[MK:M]$, \ let \ $K = \bigsqcup x_iK_1$ \ be a coset decomposition of $K$,
where \ $x_i \in K$,  $i \in {\cal I}$, and \ $x_i^{-1} x_j \in K_1$ \ if and only if $i=j$. By definition \ $|{\cal I}| = [K:K_1]$. \  
 Now consider \ $\bigcup x_iM \subseteq MK$. \ Now \ $x_i^{-1} x_j \in M$ \ implies
\ $x_i^{-1} x_j \in K\cap M = K_1$, which implies $i=j$. Hence \ 
$\bigcup x_iM = \bigsqcup x_iM \subseteq MK $. Further, if \ $g\in MK = KM$,\
then \ $g= km = x_i k_1 m$ for some $k\in K, m\in M$, and \ $k = x_ik_1$ \ for
some $i \in {\cal I}$ and some $k_1 \in K_1$ (given $k$, $i$ and $k_1$ are unique). Thus
\ $g = x_i(k_1m) \in \bigsqcup x_iM$, showing that \ $MK = \bigsqcup x_iM$ \ from which
\ $[MK:M] = |{\cal I}| = [K:K_1]$. \qed

\end{Proof}

\begin{Cor}
Let $G$ be a virtual $p$-group. Then any normal subgroup $K$ or any quotient group $H$ is also
a virtual $p$-group.

\end{Cor}

\begin{Proof}
Given a normal subgroup $K$ let $H = G/K$, or given a quotient group $H$ let \ $K$ equal the 
kernel of the projection map \ $G \twoheadrightarrow H$. Also, let $M$ be a p-subgroup of $G$ with
finite index. This gives exactly the situation of Lemma 4.6, and since \ $[G:M]$ \ is finite
the lemma implies that both $[K:K_1], \ [H:H_1]$ are finite. Clearly $K_1$ and $H_1$ are
$p$-groups since $M$ is a $p$-group. Thus $K$ and $H$ are virtual 
$p$-groups. \qed 
\end{Proof}

\begin{Thm}
 Let $G \leq A\times B $ be a virtual $p$-group. 

(a) \ If $M$ is a \syl \ of $G$, then $\overline{M}_i$ is a \syl \ of   
 $\overline{G}_i$ and $M_i$ is a \syl \ of $G_i$.

(b) \ Conversely, if $N \leq G$ has Goursat quintuple 
$(\overline{N}_1, N_1, \overline{N}_2, N_2, \theta_N)$ with
 $\overline{N}_i$  a \syl \ of   
 $\overline{G}_i$ and $N_i$  a \syl \ of $G_i$, then $N$ is a \syl \ of $G$.
\end{Thm}

\begin{Proof}
First, for both parts of the proof, notice that $G$ being a virtual $p$-group
implies, by Remark 2.3(d) and Corollary 4.7, that \ $\overline{G}_i$ and 
$G_i \approx\iota_i(G_i)  , 
\ i = 1,2$, \ are also virtual $p$-groups. 

(a)\  Let $M$ be a \syl \ of $G$, then \ $[G:M]$ \ is coprime to $p$ by 
Proposition 4.5. Applying Lemma 4.6 to Remark 2.3(d), with $\pi$ replaced by $\pi_1$, $H_1, H$ 
replaced respectively by $\overline{M}_1, \overline{G}_1$, and $K_1, K$
replaced respectively (up to isomorphism) by $M_2, G_2$, we obtain
\ $[G:M] = [\overline{G}_1 : \overline{M}_1]\cdot [G_2 : M_2]$. Hence \
$[\overline{G}_1 : \overline{M}_1] , [G_2 : M_2]$ \ are also coprime to $p$.
Since $M$ is a $p$-group so are $\overline{M}_1, M_2$, and since we have
already observed that 
$\overline{G}_1, G_2$ are virtual $p$ groups, Proposition 4.5 implies that
$\overline{M}_1, M_2$ are respectively \syls \ of $\overline{G}_1, G_2$. 
The proof for $\overline{M}_2, M_1$ is similar.

(b) \ Conversely, given the hypotheses of (b), $N$ must first of all be
a $p$-group by Proposition 4.1(c).  The steps for proving (a) above can
now all
be reversed to show that   $[G : N]$ is coprime to $p$, and hence by
Proposition 4.5 $N$ is a \syl \ of $G$. \qed

\end{Proof}

Using the correspondence between the symmetric and asymmetric versions
of Goursat's lemma (Section 2) and Theorem 3.2, the generalization of 
Theorem 4.8 to subgroups of a finite direct product is clear and we
simply state it here without proof.

\begin{Thm}
Let $G \leq A_1\times\cdots\times A_n$ be a virtual $p$-group
 and have  Goursat decomposition 
$$ Q_n(G) =   \{\overline{G}_1,\overline{G}_2, G(2|1), \theta_1,
\ldots, \overline{G}_n, G(n|1,\ldots ,n-1), \theta_{n-1}  \}  $$
as in Theorem 3.2. 

(a) If $M$ is a \syl  of $G$, then 
$\overline{M}_i$ is a \syl of 
$\overline{G}_i$ and $M(i|1,\cdots i-1)$ is a \syl of $G(i|1,\cdots i-1)$,
\ $2\leq i\leq n$.

(b)  Conversely, if $N \leq G$ has Goursat (2n-3)-tuple  
$$   \{\overline{N}_1,\overline{N}_2, N(2|1), \varphi_1,\overline{N}_3, N(3|1,2), \varphi_2 ,
\ldots, \overline{N}_n, N(n|1,\ldots ,n-1), \varphi_{n-1}  \}  $$
with each 
$\overline{N}_i$  a \syl of 
$\overline{G}_i$ and each $N(i|1,\cdots i-1)$  a \syl of $G(i|1,\cdots i-1),
\ 2\leq i\leq n$, then $N$ is a \syl of $G$.
\end{Thm}

We close this section with two questions.
\begin{Rem}
(a) \ \ One can ask whether Theorems 4.8, 4.9 hold under the weaker hypothesis
that $G$ is an ${\rm FC}_p$ group. The answer is {\it no}. Let $G$ be the
infinite cyclic subgroup of \ $\ZZ\times \ZZ_3$ \ generated by $(1,1)$, $p=3$,
and recall that any abelian group satisfies ${\rm FC}_p$. The Sylow $3$-
subgroup $M$ of $G$ equals $\{ 0\}$, hence $\overline{M}_1 = \{ 0\}$ whereas
$\overline{G}_1 = \ZZ_3$ has Sylow 3-subgroup $\ZZ_3$.

(b) \ \ In Corollary 4.7 we have seen that any normal subgroup of a virtual
$p$-group is also a virtual $p$-group. Does this hold for any subgroup?

\end{Rem}


\section{Appendix : An example illustrating the necessity of the asymmetric Goursat lemma}\label{s:Appendix}

When moving from Goursat's lemma to the asymmetric version of Goursat's lemma in Section \ref{sec:2}, we used the first isomorphism theorem to show that the required isomorphism 
\[\theta\colon \overline{G}_1/G_1\xrightarrow{\approx} \overline{G}_2/G_2\]
corresponds uniquely to a surjection
\[\theta_1\colon \overline{G}_1\twoheadrightarrow \overline{G}_2/G_2.\]
Now that we have derived
 Goursat's lemma for $n\geq 2$,
  it is tempting to use analogous
  reasoning to try to obtain a symmetric version of the lemma for $n\geq3$.  For $n=3$, such a lemma would make use of the three lattices of subgroups, each subgroup being normal in the one above it (and using the notation of Section \ref{sec:3}):
 \[ \xymatrix@!0{& \overline{G}_1\ar@{-}[dr] \ar@{-}[dl] &&&&\overline{G}_2\ar@{-}[dr] \ar@{-}[dl]&&&& \overline{G}_3\ar@{-}[dr] \ar@{-}[dl]\\
G(1|2) \ar@{-}[dr]&& G(1|3)\ar@{-}[dl]& &G(2|3)\ar@{-}[dr] && G(2|1)\ar@{-}[dl]&& G(3|1)\ar@{-}[dr] && G(3|2)\ar@{-}[dl]\\
&G(1|2,3) &&&& G(2|1,3) &&&& G(3|1,2)}\]
together with the isomorphisms
\[ \theta:\overline{G}_1/G(1|2) \to \overline{G}_2/G(2|1)\]
\[ \phi:\overline{G}_2/G(2|3) \to \overline{G}_3/G(3|2)\]
\[ \psi:\overline{G}_3/G(3|1) \to \overline{G}_1/G(1|3)\]
and  also the isomorphisms
\[ \tilde{\theta}: G(1|3)/G(1|2,3) \to G(2|3)/G(2|1,3)\]
\[ \tilde{\phi}:G(2|1)/G(2|1,3)\to G(3|1)/G(3|1,2)\]
\[ \tilde{\psi}:G(3|2)/G(3|1,2)\to G(1|2)/G(1|2,3).\]
The desired lemma would then state that this information uniquely determines the original subgroup $G$ of $A\times B\times C$.  In fact, this is {\it not} 
the case.  Indeed, we now give an example of two subspaces
of the product $A\times B\times C$ of three vector spaces, each of dimension at least $2$ (say
over $\mathbb{R}$ or $\mathbb{Q}$), which generate all of the same data as given above but nevertheless are not
the same subspace, thereby showing that a symmetric version of the lemma for $n= 3$ is
impossible. We remark that by applying forgetful functors this could also be considered
as an example in the category of abelian groups as well as the category of groups.

Choose linearly independent vectors $a_1$ and $a_2$ in $A$, $b_1$ and $b_2$ in $B$, and $c_1$ and $c_2$ in $C$. Consider the $3$-dimensional subspace
\[ V= {\operatorname{Span}}\{ (a_1, b_1, c_1), (a_2, b_2, c_2), (a_1+a_2, b_1+2b_2, c_1+3c_2)\} \]
 of $A\times B\times C$, and similarly define a second $3$-dimensional 
subspace by
\[ V' = {\operatorname{Span}}\{ (a_1, b_1, c_1+2c_2), (a_1, b_1-b_2, c_1), (a_1+2a_2, b_1, c_1)\}.\]
To compare these subspaces, we determine the twelve subspaces and nine isomorphisms presented in the subspace lattice at the beginning of this section.  For example, we find that for the subspace $V$, an element $(a,b,c)$ of the subspace 
$G(1|2)=\{ a\in A | (a, 0, c) \in V \ \text{for \ some\ } c\in C\}$ is determined by the existence of scalars $x$, $y$ and $z$ such that
\[ a= xa_1+ya_2 + z(a_1+a_2)\]
and
\[ 0=xb_1+yb_2+z(b_1+2b_2).\]
Since $b_1$ and $b_2$ are linearly independent, this leads us to the equations
\[x+z=0 \qquad \text{and} \qquad y+2z=0\]
whose solutions are ${x=y/2= -z}$.  Thus, $a= xa_1+ya_2 + z(a_1+a_2)=xa_2$ and $G(1|2)=\operatorname{Span}\{a_2\}$.
In a similar fashion, we find the following subspaces of $V$:
\[ \begin{array}{ccc}
\overline{G}_1=  \operatorname{Span}\{a_1, a_2\},&\overline{G}_2=   \operatorname{Span}\{b_1, b_2\}, &\overline{G}_3 =   \operatorname{Span}\{c_1, c_2\},\\
G(1|2)=  \operatorname{Span}\{a_2\},&G(2|3)=   \operatorname{Span}\{b_2\},&G(3|1)=   \operatorname{Span}\{c_2\}, \\
G(1|3)=  \operatorname{Span}\{ a_2\},&G(2|1)=   \operatorname{Span}\{b_2\},&G(3|2)=  \operatorname{Span}\{c_2\},  \\
G(1|2,3) =  0,&G(2| 1,3) =  0,&G(3|1,2)= 0.\\
\end{array}
\]
If we are to generalize the method of Goursat's Lemma \ref{2.1}, the isomorphism $\theta$ is determined by $\theta([a])=[b]$ where $(a,b,c)\in V$ for some $c\in C$.  Similarly, $\tilde{\theta}(a)=b$ where $(a,b,0)\in V$ (we write $\tilde{\theta}(a)$ rather than $\tilde{\theta}([a])$ since $G(1|3)/G(1|2,3)=G(1|3)/0=G(1|3)$, and likewise for $b$).  The other isomorphisms $\phi$, $\psi$, $\tilde{\phi}$ and $\tilde{\psi}$ are defined in the same manner.  The fact that this produces well-defined isomorphisms is automatically determined by the definitions of the 12 subspaces in the subspace lattice.  We wish to very specifically determine these isomorphisms for $V$, so that they can be compared to the data for $V'$.    In this case, $\theta([a_1]) = [b_1]$ since $(a_1, b_1, c_1)\in V$ and $\theta([a_2])=[b_2]=[0]$ since $(a_2, b_2, c_2)\in V$.
The other isomorphisms are determined similarly, and we obtain:
\[ \begin{array}{ccc}
\theta([a_1]) = [b_1], & \phi([b_1])= [c_1],& \psi([c_1])= [b_1],\\
\theta([a_2]) = [0], & \phi([b_2]) =  [0],& \psi([c_2]) = [0],\\
\tilde{\theta}(a_2) = b_2/2, & \tilde{\phi}(b_2) = 3c_2,& \tilde{\psi}(c_2)=-a_2.
\end{array}
\]
The determination of $\tilde{\theta}(a_2)$, for example, is given by the fact that \[(a_2, b_2/2, 0 ) = \frac{1}{2}\Big[ (a_1, b_1, c_1)+3(a_2, b_2, c_2)-(a_1+a_2, b_1+2b_2, c_1+3c_2)\Big]\] is a vector in $V$.

We leave it as an exercise for the reader to check that all twelve subspaces and six isomorphisms corresponding to the subspace $V'$ are exactly the same as those corresponding to $V$.

However, the subspaces $V$ and $V'$ are not the same.  We offer two explanations.  First, one can apply Theorem \ref{3.2} to the subspaces $V$ and $V'$.  When applying this theorem, one sees from the computations above that $\overline{G}_1$, $\overline{G}_2$, $G(2|1)$, $\theta_1$, $\overline{G}_3$ and $G(3|1,2)$ are the same for each of $V$ and $V'$.  However, the homomorphisms $\theta_2:\Gamma_2\to \overline{G}_3/G(3|1,2)$ are not the same.  The subspace $\Gamma_2$ of $A\times B$ is determined by

\[ \Gamma_2=\{ (a,b) |  (a,b,c)\in V \ (\text{resp.} V') \ \text{for\ some}\ c\in C\}\]
and the homomorphism $\theta_2$ is defined by $\theta_2(a,b)=c$, where $(a,b,c)\in V$ (resp. $V'$).
 For $V$, one finds that $(a_1, b_1)\in \Gamma_2$ with $\theta_2(a_1, b_1)=c_1$ since $(a_1, b_1, c_1)\in V$.  However, for $V'$, one finds that $(a_1, b_1)\in \Gamma_2$ but $\theta_2(a_1, b_1)=c_1+2c_2$ since $(a_1, b_1, c_1+2c_2)\in V'$.  Note that since $G(3|1,2)=0$ for both $V$ and $V'$, there is no indeterminacy so $c_1\neq c_1+2c_2$ in $\overline{G}_3$.  By Theorem \ref{3.2} we now conclude that $V\neq V'$.
 
 A more elementary explanation is given by considering a particular case of the above example in which $V$ and $V'$ are both in the Grassmann manifold $G_{6,3}$ of $3$-planes in ${\mathbb R}^6={\mathbb R}^2\oplus {\mathbb R}^2 \oplus {\mathbb R}^2$.  Let $A=B=C={\mathbb R}^2$ and let $a_1=b_1=c_1=(1,0)$ and $a_2=b_2=c_2=(0,1)$.  Then the subspace $V$ is the span of \[\{ (1,0,1,0,1,0), (0,1,0,1,0,1), (0,1,0,2,0,3)\}\] in ${\mathbb R}^6$, while $V'$ is the span of \[\{ (1,0,1,0,1,2), (1,0,1,-1,1,0), (1,2,1,0,1,0)\}.\]  The matrix
 \[ \left[ \begin{array}{cccccc}
 1 & 0 & 1 & 0 & 1 & 0\\
 0 & 1 & 0 & 1 & 0 & 1\\
 0 & 1 & 0 & 2 & 0 & 3\\
 1 & 0 & 1& 0 & 1 & 2
 \end{array}
 \right]\]
 has rank 4, showing that the vector $(1, 0, 1, 0, 1, 2)$ of $V'$ is not in $V$.  Thus $V\neq V'$ in this case.

This demonstrates that the symmetric version of Goursat's lemma does not hold.

\begin{Rem} Upon trying to recover the symmetric version of Goursat's lemma by applying the first isomorphism theorem to Theorem \ref{3.2}, one sees almost immediately that the trouble stems from the fact that $G(1|2,3)$ is not necessarily equal to the intersection of $G(1|2)$ and $G(1|3)$, and similarly for $G(2|1,3)$ and $G(2|1,3)$.  It seems very likely that a symmetric version of Goursat's lemma is available with this additional hypothesis.
\end{Rem}

\providecommand{\bysame}{\leavevmode\hbox to3em{\hrulefill}\thinspace}
\providecommand{\MR}{\relax\ifhmode\unskip\space\fi MR }
\providecommand{\MRhref}[2]{%
  \href{http://www.ams.org/mathscinet-getitem?mr=#1}{#2}
}
\providecommand{\href}[2]{#2}



\noindent
\author{Kristine Bauer\\
{\small {Department of Mathematics and Statistics} \\
\small{University of Calgary} \\
\small {Calgary, Alberta, Canada T2N 1N4}\\
\small {e-mail: \url{bauerk@ucalgary.ca}}\\

\noindent
\author{Debasis Sen} \\
{\small {Department of Mathematics and Statistics} \\
\small{Indian Institute of Technology} \\
\small {Kanpur, India}\\
\small {e-mail: \url{debasis@iitk.ac.in}}\\

\noindent
\author{Peter Zvengrowski} \\
{\small {Department of Mathematics and Statistics} \\
\small{University of Calgary} \\
\small {Calgary, Alberta, Canada T2N 1N4}\\
\small { e-mail: \url{zvengrow@ucalgary.ca}}

\end{document}